\pdfoutput=1
\RequirePackage{ifpdf}
\ifpdf 
\documentclass[pdftex]{sigma}
\else
\documentclass{sigma}
\fi

\def\Z{\mathbb Z}
\def\Q{\mathbb Q}

\def\C{\mathbb C}
\def\P{\mathbb P}

\def\M{{\mathcal M}}

\def\X{\mathcal X}

\def\L{\mathcal L}

\def\<{\langle}
\def\>{\rangle}

\def\a{\alpha}
\def\b{\beta}

\def\embd{\hookrightarrow}
\def\Aut{\operatorname{Aut}}

\numberwithin{equation}{section}
\newtheorem{thm}{Theorem}[section]
\newtheorem{prop}[thm]{Proposition}
\newtheorem{lem}[thm]{Lemma}
\newtheorem{cor}[thm]{Corollary}

\theoremstyle{definition}

\newtheorem{rem}[thm]{Remark}

\newcommand{\im}{\operatorname{Im}}

\newcommand{\supp}{\operatorname{supp}}

\newcommand{\sltwoz}{\operatorname{SL}_2(\Z)}
\newcommand{\sltwoc}{\operatorname{SL}_2(\C)}

\newcommand{\gltwoc}{\operatorname{GL}_2(\C)}
\newcommand{\sptwoz}{\operatorname{Sp}_4(\Z)}

\newcommand{\spgz}{\operatorname{Sp}_{2g}(\Z)}

\def\x{\mathbf{x}}
\def\y{\mathbf{y}}
\def\p{\mathfrak p}
\def\wM{\mathbb{WP}^3_{(2,4,6,10)}}

\begin{document}
\allowdisplaybreaks

\newcommand{\arXivNumber}{1607.08294}

\renewcommand{\thefootnote}{}

\renewcommand{\PaperNumber}{089}

\FirstPageHeading

\ShortArticleName{A Universal Genus-Two Curve from Siegel Modular Forms}

\ArticleName{A Universal Genus-Two Curve\\ from Siegel Modular Forms\footnote{This paper is a~contribution to the Special Issue on Modular Forms and String Theory in honor of Noriko Yui. The full collection is available at \href{http://www.emis.de/journals/SIGMA/modular-forms.html}{http://www.emis.de/journals/SIGMA/modular-forms.html}}}

\Author{Andreas MALMENDIER~$^\dag$ and Tony SHASKA~$^\ddag$}

\AuthorNameForHeading{A.~Malmendier and T.~Shaska}

\Address{$^\dag$~Department of Mathematics and Statistics, Utah State University, Logan, UT 84322, USA}
\EmailD{\href{mailto:andreas.malmendier@usu.edu}{andreas.malmendier@usu.edu}}

\Address{$^\ddag$~Department of Mathematics and Statistics, Oakland University, Rochester, MI 48309, USA}
\EmailD{\href{mailto:shaska@oakland.edu}{shaska@oakland.edu}}

\ArticleDates{Received July 18, 2017, in f\/inal form November 25, 2017; Published online November 30, 2017}

\Abstract{Let $\p $ be any point in the moduli space of genus-two curves $\M_2$ and $K$ its f\/ield of moduli. We provide a universal equation of a genus-two curve $\mathcal C_{\a, \b}$ def\/ined over $K(\alpha, \beta)$, corresponding to $\p$, where $\alpha $ and $\beta$ satisfy a quadratic $\alpha^2+ b \beta^2= c$ such that $b$ and $c$ are given in terms of ratios of Siegel modular forms. The curve $\mathcal C_{\a, \b}$ is def\/ined over the f\/ield of moduli $K$ if and only if the quadratic has a $K$-rational point $(\a, \b)$. We discover some interesting symmetries of the Weierstrass equation of $\mathcal C_{\a, \b}$.
This extends previous work of Mestre and others.}


\Classification{14H10; 14H45}

\renewcommand{\thefootnote}{\arabic{footnote}}
\setcounter{footnote}{0}

\section{Introduction}

Let $\M_2$ be the moduli space of genus-two curves. It is the coarse moduli space for smooth, complete, connected curves of genus two over $\C$.
Let $\p \in \M_2 (K)$, where $K$ is the f\/ield of def\/inition of~$\p$. Construction a genus 2 curve $C$ corresponding to $\p$ is interesting from many points of view. Mestre \cite{Me} has shown how to construct equations for genus-two curves with automorphism group of order two and def\/ined over $\Q$. Mestre's work has recently received new attention from researchers in experimental number theory. For instance, in~\cite{MR3540958} a database of geometric and arithmetic invariants of genus-two curves def\/ined over $\Q$ of small discriminant. In \cite{data}, the authors count the points in $\M_2 (\Q)$ according to their moduli height and create a database of genus-two curves from the moduli points in $\M_2 (\Q)$. In creating the database the main problem was that of constructing an equation for obstruction moduli points. This paper provides an equation over a minimal f\/ield of def\/inition for \emph{any} point $\p \in \M_2$. Our work is therefore complimentary to the problem of f\/inding an ef\/f\/icient construction for genus-two curves over f\/inite f\/ields with a prescribed number of rational points and the associated complexity analysis in \cite{MR3349314, MR2899960}. Our equation for a genus-two curve is universal in the sense that it works for every moduli point given in terms of Igusa invariants or Siegel modular forms. It does not rely on special CM values for Siegel modular functions where the associated abelian surface has extra endomorphisms or the special invariants that can be used in theses cases (cf.~\cite{MR3435723}).

The natural question is if there exists a universal curve for the genus-two curve given in terms of a generic moduli point $\p\in \M_2$. In other words, given an af\/f\/ine moduli point $\p = (x, y, z)$, where $x$, $y$, $z$ are transcendentals, can we construct a curve corresponding to $\p$? The answer is
negative in the strict def\/inition of ``universal curve''; see \cite[p.~39]{mod-curves} for details. As we will show, there is a satisfactory answer in the sense that our ``universal equation" applies to every moduli point $\p \in \M_2$. However, the equation is often def\/ined only over a quadratic extension of the f\/ield of moduli.

We focus mainly on constructing a genus-two curve $\mathcal C$ for any given point $\p=(\x_1, \x_2, \x_3)\in \M_2$, def\/ined over a minimal f\/ield of def\/inition, where $\x_1$, $\x_2$, $\x_3$ are ratios of modular forms as def\/ined by Igusa in~\cite{MR0141643}.
Our main result is as follows: For every point $\p \in \M_2$ such that $\p \in \M_2 (K)$, where $K$ is the f\/ield of moduli, there is a genus-two curve $\mathcal{C}_{( \alpha, \beta)}$ given by
\begin{gather*}
 \mathcal{C}_{( \alpha, \beta)}\colon \ y^2= \sum_{i=0}^6 a_{i}( \alpha, \beta) x^i ,
\end{gather*}
corresponding to $\p$ with coef\/f\/icients given by equation~\eqref{SexticPolynomial}. This curve is def\/ined over the f\/ield of moduli $K$ if and only if there exists a $K$-rational solution $(\a, \b)$ to the quadratic
\begin{gather*} \a^2 + b \cdot \b^2 = c\end{gather*}
where $b$ and $c$ are given in terms of the moduli point $\p$. There are some interesting properties of the coef\/f\/icients def\/ining $\mathcal{C}_{( \alpha, \beta)}$ which seem to be particular to this model and not noticed before.

It must be noticed that this equation is universal in the sense that it works for every moduli point $[J_2: J_4: J_6: J_{10}]$ given in terms of the
Igusa invariants $J_2$, $J_4$, $J_6$, $J_{10}$. The equation is def\/ined at worst over a quadratic extension of the f\/ield of moduli $K$. If the equation over the f\/ield of moduli is needed, then we must search locally for a rational point in the above quadratic when evaluated at the given $\p$. In the process we discover some
interesting absolute invariants (cf. equation~\eqref{OurCoordinates}) which as far as we are aware have not been used before.

The paper is organized as follows: In Section~\ref{section-two} we give a brief summary of Siegel modular forms, classical invariants of binary sextics and the relations among them. While this material can be found in many places in the literature, there is plenty of confusion on the labeling
and normalization of such invariants and relations among them.
We also introduce a set of absolute invariants that is well-suited for the construction of a universal sextic.

In Section~\ref{universal-curve} we construct the equation of the genus-two curve by determining the Clebsch conic and the cubic. We diagonalize the corresponding conic and discover a new set of invariants which make the equation of this conic short and elegant. The diagonalized conic can be quickly determined from the invariants of the curve. The intersection of this conic and the cubic gives the equation of the genus-two curve. This equation shows some interesting symmetries of the coef\/f\/icients, which to the knowledge of the authors have never been discovered before. When this universal equation is restricted to loci of curves with automorphisms or the Clebsch invariant $D=0$ (not covered by Mestre's approach) it shows that the f\/ield of moduli is a f\/ield of def\/inition, results which agree with previous results of other authors.

\section{Preliminaries}\label{section-two}

\subsection{The Siegel modular three-fold} \label{SiegelThreefold}
The Siegel three-fold is a quasi-projective variety of dimension $3$ obtained from the Siegel upper
half-plane of degree two which by def\/inition is the set of two-by-two symmetric matrices over $\C$ whose imaginary part is positive def\/inite, i.e.,
\begin{gather*}
 \mathbb{H}_2 = \left. \left\lbrace \underline{\tau} = \left( \begin{matrix} \tau_1 & z \\ z & \tau_2\end{matrix} \right) \right|
 \tau_1, \tau_2, z \in \C , \,\im{(\tau_1)} \im{(\tau_2}) > \im{(z)}^2 , \,\im{(\tau_2)} > 0 \right\rbrace ,
\end{gather*}
quotiented out by the action of the modular transformations $\Gamma_2:=\sptwoz$, i.e.,
\begin{gather*}
 \mathcal{A}_2 = \mathbb{H}_2 / \Gamma_2 .
\end{gather*}
Each $\underline{\tau} \in \mathbb{H}_2$ determines a principally polarized complex abelian surface
\begin{gather*} \mathbf{A}_{\underline{ \tau}} = \C^2 / \big\langle \Z^2 \oplus \underline{\tau} \Z^2\big\rangle\end{gather*}
 with period matrix $(\underline{\tau}, \mathbb{I}_2) \in \operatorname{Mat}(2, 4;\C)$.
Two abelian surfaces $\mathbf{A}_{\underline{ \tau}}$ and $\mathbf{A}_{\underline{ \tau}'}$ are isomorphic if and only if there is a symplectic matrix
\begin{gather*}
M= \left(\begin{matrix} A & B \\ C & D \end{matrix} \right) \in \Gamma_2,
\end{gather*}
such that $\underline{\tau}' = M (\underline{\tau}):=(A\underline{\tau}+B)(C\underline{\tau}+D)^{-1}$.
It follows that the Siegel three-fold $\mathcal{A}_2$ is also the set of isomorphism classes of principally polarized abelian surfaces.
The sets of abelian surfaces that have the same endomorphism ring form sub-varieties of $\mathcal{A}_2$.
The endomorphism ring of principally polarized abelian surface tensored with $\mathbb{Q}$ is either a quartic CM f\/ield, an indef\/inite
quaternion algebra, a real quadratic f\/ield or in the generic case $\mathbb{Q}$. Irreducible components of the corresponding
subsets in $\mathcal{A}_2$ have dimensions $0$, $1$, $2$ and are known as CM points, Shimura curves and
Humbert surfaces, respectively.

The Humbert surface $H_{\Delta}$ with invariant $\Delta$ is the space of principally polarized abelian surfaces admitting a symmetric endomorphism
with discriminant $\Delta$. It turns out that $\Delta$ is a~positive integer $\equiv 0, 1\mod 4$. In fact, $H_{\Delta}$ is the image inside $\mathcal{A}_2$
under the projection of the rational divisor associated to the equation
\begin{gather*}
 a \tau_1 + b z + c \tau_3 + d \big(z^2 -\tau_1 \tau_2\big) + e = 0 ,
\end{gather*}
with integers $a$, $b$, $c$, $d$, $e$ satisfying $\Delta=b^2-4 a c-4 d e$ and $\underline{\tau}
= \bigl(\begin{smallmatrix}
\tau_1&z\\ z&\tau_2
\end{smallmatrix} \bigr) \in \mathbb{H}_2$.
For example, inside of $\mathcal{A}_2$ sit the Humbert surfaces $H_1$ and $H_4$ that are def\/ined as the images under the projection of the rational divisor associated to $z=0$ and $\tau_1=\tau_2$, respectively. In fact, the singular locus of $\mathcal{A}_2$ has~$H_1$ and~$H_4$ as its two connected components.
As analytic spaces, the surfaces~$H_1$ and~$H_4$ are each isomorphic to the Hilbert modular surface
\begin{gather*}
 \big( (\sltwoz \times \sltwoz) \rtimes \Z_2 \big) \backslash \big( \mathbb{H} \times \mathbb{H} \big) .
\end{gather*}
For a more detailed introduction to Siegel modular form, Humbert surfaces, and the Satake compactif\/ication of the Siegel modular threefold we refer to
Freitag's book \cite{MR871067}.

\subsection{Siegel modular forms}\label{Siegel_modular_forms}

In general, we can def\/ine the Eisenstein series $\psi_{2k}$ of degree $g$ and weight $2k$ (where we assume $2k>g+1$ for convergence) by setting
\begin{gather*}
 \psi_{2k}(\underline{\tau}) = \sum_{(C,D)} \det(C\cdot\underline{\tau}+D)^{-2k} ,
\end{gather*}
where the sum runs over non-associated bottom rows $(C,D)$ of elements in $\spgz$ where non-associated means with respect to the multiplication by $\mathrm{GL}_g(\Z)$. In the following, we will always assume $g=2$ in the def\/inition of $\psi_{2k}$.
Using Igusa's def\/inition \cite[Section~8, p.~195]{MR0141643} we def\/ine a cusp form of weight $10$ by
\begin{gather*}
 \chi_{10}(\underline{\tau}) = - \frac{43867}{2^{12} 3^5 5^2 7 \cdot 53} \left(\psi_4(\underline{\tau}) \psi_6(\underline{\tau}) - \psi_{10}(\underline{\tau})\right) .
\end{gather*}
Based on Igusa's def\/inition \cite[Section~8, p.~195]{MR0141643} and the work in \cite{MR3366121} we def\/ine a second cusp form $\chi_{12}$ of weight $12$ by
\begin{gather*}
 \chi_{12}(\underline{\tau}) = \frac{131\cdot 593}{2^{13} 3^7 5^3 7^2 337}\big(3^2 7^2 \psi_4^3(\underline{\tau})+2\cdot 5^3 \psi_6^2(\underline{\tau})-691 \psi_{12}(\underline{\tau})\big) .
\end{gather*}
Moreover, Igusa proved \cite{MR0229643, MR527830} that the ring of Siegel modular forms is generated by
$\psi_4$, $\psi_6$, $\chi_{10}$, $\chi_{12}$ and by one more cusp form $\chi_{35}$ of odd weight $35$
whose square is the following polynomial \cite[p.~849]{MR0229643} in the even generators
\begin{gather*}
\chi_{35}^2  = \frac{1}{2^{12} 3^9} \chi_{10} \big(
2^{24} 3^{15} \chi_{12}^5 - 2^{13} 3^9 \psi_4^3 \chi_{12}^4 - 2^{13} 3^9 \psi_6^2 \chi_{12}^4 + 3^3 \psi_4^6 \chi_{12}^3 - 2\cdot 3^3 \psi_4^3 \psi_6^2 \chi_{12}^3 \\
\hphantom{\chi_{35}^2  =}{} - 2^{14} 3^8 \psi_4^2 \psi_6 \chi_{10} \chi_{12}^3 -2^{23} 3^{12} 5^2 \psi_4 \chi_{10}^2 \chi_{12}^3 + 3^3 \psi_6^4 \chi_{12}^3+ 2^{11} 3^6 37 \psi_4^4 \chi_{10}^2 \chi_{12}^2\\
\hphantom{\chi_{35}^2  =}{}+2^{11} 3^6 5\cdot 7 \psi_4 \psi_6^2 \chi_{10}^2 \chi_{12}^2 -2^{23} 3^9 5^3 \psi_6 \chi_{10}^3 \chi_{12}^2- 3^2 \psi_4^7 \chi_{10}^2 \chi_{12} + 2 \cdot 3^2 \psi_4^4 \psi_6^2 \chi_{10}^2 \chi_{12} \\
\hphantom{\chi_{35}^2  =}{} + 2^{11} 3^5 5 \cdot 19 \psi_4^3 \psi_6 \chi_{10}^3 \chi_{12}  + 2^{20} 3^8 5^3 11 \psi_4^2 \chi_{10}^4 \chi_{12} - 3^2 \psi_4 \psi_6^4 \chi_{10}^2 \chi_{12} + 2^{11} 3^5 5^2 \psi_6^3 \chi_{10}^3 \chi_{12} \\
\hphantom{\chi_{35}^2  =}{} - 2 \psi_4^6 \psi_6 \chi_{10}^3 - 2^{12} 3^4 \psi_4^5 \chi_{10}^4 + 2^2 \psi_4^3 \psi_6^3 \chi_{10}^3 + 2^{12} 3^4 5^2 \psi_4^2 \psi_6^2 \chi_{10}^4 + 2^{21} 3^7 5^4 \psi_4 \psi_6 \chi_{10}^5 \\
\hphantom{\chi_{35}^2  =}{}  - 2 \psi_6^5 \chi_{10}^3 + 2^{32} 3^9 5^5 \chi_{10}^6 \big) .
\end{gather*}
Hence, the expression $Q:= 2^{12} 3^9 \chi_{35}^2 /\chi_{10}$ is a polynomial of degree $60$ in the even generators. Igusa also proved that each Siegel modular form (with trivial character) of odd weight is divisible by the form $\chi_{35}$. The following fact is known\cite{MR1438983}:

\begin{prop}The vanishing divisor of $Q$ is the Humbert surface $H_4$, i.e., a period point~$\underline{\tau}$ is equivalent to a point with $\tau_1=\tau_2$ if and only if $Q=0$. Accordingly, the vanishing divisor of~$\chi_{35}$ is the formal sum $H_1 + H_4$ of Humbert surfaces, that constitutes the singular locus of~$\mathcal{A}_2$.
\end{prop}

In accordance with Igusa \cite[Theorem~3]{MR0141643} we also introduce the following ratios of Siegel modular forms
\begin{gather*}
 \x_1 = \dfrac{\psi_4 \chi_{10}^2}{\chi_{12}^2} ,\qquad
 \x_2 = \dfrac{\psi_6 \chi_{10}^3}{\chi_{12}^3} ,\qquad
 \x_3 = \dfrac{\chi_{10}^6}{\chi_{12}^5} ,
\end{gather*}
as well as
\begin{gather}\label{ratio_of_Siegel_forms2}
 \y_1 = \dfrac{\x_1^3}{\x_3} = \dfrac{\psi_4^3}{\chi_{12}} ,\qquad
 \y_2 = \dfrac{\x_2^2}{\x_3} = \dfrac{\psi_6^2}{\chi_{12}} ,\qquad
 \y_3 = \dfrac{\x_1^2 \x_2}{\x_3} = \dfrac{\psi_4^2 \psi_6 \chi_{10}}{\chi_{12}} ,
\end{gather}
where we have suppressed the dependence of each Siegel modular form on $\underline{\tau}$. These ratios have the following asymptotic expansion for $z\to 0$ \cite[pp.~180--182]{MR0141643} in terms of ordinary Eisenstein series $E_4$ and $E_6$ and the Dedekind $\eta$-function
\begin{gather*}
 \x_1 = E_4(\tau_1) E_4(\tau_2) (\pi z)^{4} + O\big(z^{5}\big) , \nonumber\\
 \x_2 = E_6(\tau_1) E_6(\tau_2) (\pi z)^{6} + O\big(z^{7}\big) , \nonumber\\
 \x_3 = \eta^{24}(\tau_1) \eta^{24}(\tau_2) (\pi z)^{12} + O\big(z^{13}\big) ,\label{asymptotics2}
\end{gather*}
and
\begin{gather}
 \y_1 = j(\tau_1) j(\tau_2) + O\big(z^2\big) ,\nonumber \\
 \y_2 = (1728 - j(\tau_1)) (1728 -j(\tau_2)) + O\big(z^2\big) , \nonumber\\
 \y_3 = \dfrac{E_4^2(\tau_1) E_4^2(\tau_2) E_6(\tau_1) E_6(\tau_2)}{\eta^{24}(\tau_1) \eta^{24}(\tau_2)} (\pi z)^2 + O\big(z^3\big) ,\label{asymptotics3}
\end{gather}
where we have set
\begin{gather*}
 j(\tau_j)  = \dfrac{1728 E_4^3(\tau_j)}{E_4^3(\tau_j)-E_6^2(\tau_j)} =\dfrac{E_4^3(\tau_j)}{\eta^{24}(\tau_j)} , \\
 1728 - j(\tau_j)  = \dfrac{1728 E_6^2(\tau_j)}{E_4^3(\tau_j)-E_6^2(\tau_j)} =\dfrac{E_6^2(\tau_j)}{\eta^{24}(\tau_j)} .
\end{gather*}

The following fact follows from the above asymptotic analysis \cite{MR1438983}:

\begin{prop}\label{Prop_chi10}
The modulus point $\underline{\tau}$ is equivalent to a point with $z=0$ or $[\underline{\tau}] \in H_1 \subset \mathcal{A}_2$ such that
the principally polarized abelian surface is a product of two elliptic curves $\mathbf{A}_{\underline{\tau}}=E_{\tau_1} \times E_{\tau_2}$
if and only if $\chi_{10}(\underline{\tau}) =0$. The elliptic modular parameters are determined by equation~\eqref{asymptotics3}.
\end{prop}

\subsection{Igusa invariants} \label{moduli_curves_genus2}
Suppose that $\mathcal{C}$ is an irreducible projective non-singular curve. If the self-intersection is $\mathcal{C}\cdot \mathcal{C}=2$ then $\mathcal{C}$ is a curve of genus two. For every curve $\mathcal{C}$ of genus two there exists a unique pair $(\operatorname{Jac}(\mathcal{C}),j_\mathcal{C})$ where $\operatorname{Jac}(\mathcal{C})$ is an abelian surface, called the Jacobian variety of the curve $\mathcal{C}$, and $j_\mathcal{C}\colon \mathcal{C} \to \operatorname{Jac}(\mathcal{C})$ is an embedding. One can always regain $\mathcal{C}$ from the pair $(\operatorname{Jac}(\mathcal{C}),\mathcal{P})$ where
$\mathcal{P}=[\mathcal{C}]$ is the class of $\mathcal{C}$ in the N\'eron--Severi group $\mathrm{NS}(\operatorname{Jac}(\mathcal{C}))$. Thus, if $\mathcal{C}$ is a genus-two curve, then $\operatorname{Jac}(\mathcal{C})$ is a principally polarized abelian surface with principal polarization $\mathcal{P}=[\mathcal{C}]$, and the map sending a curve $\mathcal{C}$ to its Jacobian variety $\operatorname{Jac}(\mathcal{C})$ is injective. In this way, the variety of moduli of curves of genus two is also the moduli space of their Jacobian varieties with canonical polarization.

We write the equation def\/ining a genus-two curve $\mathcal{C}$ by a degree-six polynomial or sextic in the form
\begin{gather}\label{genus_two_curve}
 \mathcal{C}: y^2 = f(x) = a_0 \prod_{i=1}^6 (x-\a_i) = \sum_{i=0}^6 a_i x^{i} .
\end{gather}
The roots $\{\a_i \}_{i=1}^6$ of the sextic are the six ramif\/ication points of the map $\mathcal{C} \to \mathbb{P}^1$. Their pre-images on $\mathcal{C}$ are the six Weierstrass points. The isomorphism class of $f$ consists of all equivalent
sextics where two sextics are considered equivalent if there is a linear transformation in $\gltwoc$ which takes the
set of roots to the roots of the other.

The ring of invariants of binary sextics is generated by the Igusa invariants $(J_2 , J_4 , J_6 , J_{10})$ as def\/ined in \cite[equation~(9)]{vishy}, which are the same invariants as the ones denoted by $( A' , B' , C' , D')$ in \cite[p.~319]{Me} and also the same invariants as $(A, B, C, D)$ in \cite[p.~176]{MR0141643}.
 For expressions of such invariants in terms of the coef\/f\/icients $a_0, \dots, a_6$ of the binary sextic, or $J_k \in \Z[a_0,\dots,a_6]$ for $k\in \{2, 4, 6, 10\}$; see \cite[equation~(11)]{vishy} and in terms of thetanulls see \cite{satake}.
One can then ask what the Igusa invariants of a genus-two curve $\mathcal{C}$ def\/ined by a sextic curve $f$ are in terms of $\underline{\tau}$
such that $(\underline{\tau}, \mathbb{I}_2) \in \mathrm{Mat}(2, 4;\mathbb{C})$ is the period matrix of the principally polarized abelian surface
$\mathbf{A}_{\underline{\tau}}=\operatorname{Jac}(\mathcal{C})$. Based on the asymptotic behavior in Equations \eqref{asymptotics2} and \eqref{asymptotics3}, Igusa \cite[p.~848]{MR0229643} proved that the relations are as follows
\begin{gather*}
 J_2  = -2^3 \cdot 3 \dfrac{\chi_{12}(\underline{\tau})}{\chi_{10}(\underline{\tau})} , \\
 J_4  = \phantom{-} 2^2 \psi_4(\underline{\tau}) ,\\
 J_6  = -\frac{2^3}{3} \psi_6(\underline{\tau}) - 2^5 \dfrac{\psi_4(\underline{\tau}) \chi_{12}(\underline{\tau})}{\chi_{10}(\underline{\tau})} ,\\
 J_{10} = -2^{14} \chi_{10}(\underline{\tau}) .\label{invariants}
\end{gather*}

Thus, the invariants of a sextic def\/ine a point in a weighted projective space $[J_2 : J_4 : J_6 : J_{10}] \in \mathbb{WP}^3_{(2,4,6,10)}$
that equals
\begin{gather*}\label{IgusaClebschProjective}
 \big[ 2^3 3 (3\chi_{12}) : 2^2 3^2 \psi_4 \chi_{10}^2 : 2^3 3^2 \big(4 \psi_4 (3\chi_{12})+ \psi_6 \chi_{10} \big) \chi_{10}^2: 2^2 \chi_{10}^6 \big] .
\end{gather*}

Torelli's theorem states that the map sending a curve $\mathcal{C}$ to its Jacobian variety $\operatorname{Jac}(\mathcal{C})$ induces a birational map
from the moduli space $\mathcal{M}_2$ of genus-two curves to the complement of the Humbert surface $H_1$ in $\mathcal{A}_2$, i.e.,
$\mathcal{A}_2 - \supp{(\chi_{10})}_0$. In other words, points in the projective variety $\operatorname{Proj} \C [J_2, J_4, J_6, J_{10}]$ which
are not on $J_{10}=0$ are in one-to-one correspondence with isomorphism classes of regular sextics \cite{MR0141643}.

Often the \emph{Clebsch invariants} $(A, B, C, D)$ of a sextic are used instead. They are def\/ined in terms of the transvectants of the binary sextics; see~\cite{MR0485106} for details.
The invariants $(A,B,C,D)$ are polynomial expressions in the Igusa invariants $(J_2, J_4, J_6, J_{10})$ with rational coef\/f\/icients:
\begin{gather}
 A = -\frac{1}{2^3 3 \cdot 5} J_2 , \nonumber\\
 B = \frac{1}{2^3 3^3 5^4} \big( J_2^2 + 20 J_4 \big) ,\nonumber\\
 C =  - \frac{1}{2^5 3^5 5^6} \big( J_2^3 + 80 J_2 J_4 - 600 J_6\big) ,\label{iClebsch_invariants}\\
 D =  - \frac{1}{2^8 3^9 5^{10}} \big( 9 J_2^5+700 J_2^3 J_4-3600 J_2^2 J_6 -12400 J_2 J_4^2+48000 J_4 J_6+10800000 J_{10}\big) .\nonumber
\end{gather}
For formulas giving relations between all these sets of invariants see \cite{data}.

\subsubsection{Absolute invariants}
Dividing any $\sltwoc$ invariant by another one of the same degree gives an invariant under $\gltwoc$ action. The term \textit{absolute invariants} is used f\/irst by Igusa \cite{MR0114819} for $\gltwoc$ invariants.
It was the main result of \cite[{Theorem 3}]{MR0141643} that
\begin{gather*}
\x_1=144 \frac {J_4} {J_2^2}, \qquad \x_2=- 1728 \frac {J_2J_4-3J_6} {J_2^3}, \qquad \x_3 =486 \frac {J_{10}} {J_2^5},
\end{gather*}
for $J_2 \neq 0$. We use $\x_1$, $\x_2$, $\x_3$ to write the point $[J_2 : J_4 : J_6 : J_{10}] \in \mathbb{WP}^3_{(2,4,6,10)}$ as
\begin{gather*}
\left\lbrack 1 : \frac{1}{2^4 3^2} \x_1 : \frac{1}{2^6 3^4} \x_2 + \frac{1}{2^4 3^3} \x_1 : \frac{1}{2 \cdot 3^5} \x_3\right\rbrack .
\end{gather*}
Since the invariants $J_4$, $J_6$, $J_{10}$ vanish simultaneously for sextics with triple roots all such curves are mapped to $[1:0:0:0] \in \mathbb{WP}^3_{(2,4,6,10)}$ with uniformizing af\/f\/ine coordinates~$\x_1$,~$\x_2$,~$\x_3$ around it. Blowing up this point gives a variety that
parameterizes genus-two curves with $J_2 \not = 0$ and their degenerations. In the blow-up space we have to introduce additional coordinates that are obtained as ratios of $\x_1$, $\x_2$, $\x_3$ and have weight zero. Those are precisely the coordinates~$\y_1$,~$\y_2$,~$\y_3$ already introduced in equation~\eqref{ratio_of_Siegel_forms2}. It turns out that the coordinate ring of the blown-up space is $\C[\x_1, \x_2, \x_3, \y_1, \y_2, \y_3]$.

We introduce the three absolute invariants
\begin{gather}
 \rho =- \frac { 4 \big( 9 J_2^{2}-320 J_4 \big) \big( J_2^{2}+20 J_4 \big) ^{2}}{ \big( 3 J_2^{3}+140 J_2 J_4-800 J_6 \big) ^{2}} ,\nonumber\\
 \sigma  = -\frac { 48 \big( J_2^{2}+20 J_4 \big) ^{2} }{ \big( 3 J_2^{3}+140 J_2 J_4-800 J_6 \big) ^{3}}\nonumber \\
\hphantom{\sigma  =}{} \times \big( 9 J_2^{5}-700 J_2^{3} J_4+2400 J_2^{2} J_6-262400 J_2 J_4^{2}+768000 J_4 J_6
 +172800000 J_{10} \big),\nonumber\\
 \kappa =\frac {2 \big( 27 J_2^{4}+2380 J_2^2 J_4-12000 J_2 J_6+12800 J_4^2 \big) \big( J_2^{2}+20 J_4 \big) }{ \big( 3 J_2^{3}+140 J_2 J_4-800 J_6 \big) ^{2}} .\label{OurCoordinates}
\end{gather}
It follows:
\begin{lem} \label{lem-NewInv}
For invariants $(\rho, \sigma, \kappa)$ given by equation~\eqref{OurCoordinates} such that $\rho$ and $\kappa$ do not vanish simultaneously,
a point $[J_2 : J_4 : J_6 : J_{10}]$ in $\mathbb{WP}^3_{(2,4,6,10)}$ is given by
\begin{gather}
J_2 = 8(\kappa-\rho), \qquad J_4 = \frac{9}{5}(\kappa-\rho)^2+45\rho, \nonumber\\
J_6= \frac{111}{25}(\kappa-\rho)^3-30(\kappa-\rho)^2+63\rho(\kappa-\rho)-270\rho ,\nonumber\\
J_{10}  = \frac{6}{3125} (\kappa-\rho)^5+\frac{4}{15} (\kappa-\rho)^4+\frac{46}{75} \rho (\kappa-\rho)^3\nonumber\\
\hphantom{J_{10}  =}{} +\left(-\frac{1}{6}\sigma+\frac{42}{5}\rho\right) (\kappa-\rho)^2+12 \rho^2 (\kappa-\rho)+\frac{3}{2} \rho (36\rho-\sigma).\label{eq-21}
\end{gather}
In particular, for $J_2 \not = 0$ we have $\Q (\x_1, \x_2, \x_3 ) = \Q (\rho, \sigma, \kappa)$.
\end{lem}
\begin{proof} The proof is computational. We express $\rho, \sigma, \kappa$ as rational functions of $\x_1$, $\x_2$, $\x_3$ and vice versa
over $\mathbb{Q}$. The condition that $\rho$ and $\kappa$ do not vanish simultaneously is based on the fact that
$J_2$, $J_4$, $J_6$, $J_{10}$ must not vanish simultaneously.
\end{proof}

\begin{rem}
Consider the image of $[J_2 : J_4 : J_6 : J_{10}]$ in $\mathbb{WP}^3_{(2,4,6,10)}$ under the morphism $\mathbb{WP}^3_{(2,4,6,10)} \to \mathbb{P}^5$ given by
\begin{gather}\label{embedding}
 \big[486 J_4 J_6: 486 J_{10}: -1728 (J_2 J_4-3 J_6) J_2^2: 144 J_2^3 J_4: 20736 J_2 J_4^2: J_2^5\big],
\end{gather}
which is a linear transformation of the usual morphism to $ \mathbb{P}^5$ given by
\begin{gather*}
[J_2 : J_4 : J_6 : J_{10}] \mapsto \big[J_{10}: J_4J_6: J_2^2 J_6: J_2^3 J_4: J_2 J_4^2: J_2^5\big].
\end{gather*}
For $J_2 \not =0$, points in equation~\eqref{embedding} equal
\begin{gather*}
 \left[\frac{1}{1536}\x_1(\x_2+12\x_1): \x_3: \x_2: \x_1: \x_1^2: 1\right].
\end{gather*}
The invariants $\x_1$, $\x_2$, $\x_3$ are not def\/ined for $J_2 =0$, but $\rho$, $\sigma$, $\kappa$ remain well-def\/ined if $\rho=\kappa\not = 0$. In this case we have
\begin{gather*}
 J^{(0)}_2 = 0, \qquad J^{(0)}_4 = 45\rho, \qquad J^{(0)}_6 = -270\rho, \qquad J^{(0)}_{10} = \frac{3}{2}\rho (36\rho-\sigma) ,
\end{gather*}
and the invariants $\rho$ and $\sigma$ with
\begin{gather*}
 \rho = \kappa = \frac{4}{5} \frac{J_4^3}{J_6^2}, \qquad \sigma = \frac{144}{5} \frac{J_4^3}{J_6^2} + 6480 \frac{J_4^3}{J_6^2} \frac{J_{10}}{J_4 J_6},
 \end{gather*}
determine genus-two curves with $J_2=0$, $J_4 \cdot J_{6}\neq 0$ up to isomorphism. In addition to $J_2=0$, we have $J_{10}=0$
if and only if $\sigma=36\rho$. Using $\epsilon=(\kappa-\rho)$ in equation~\eqref{eq-21}, one checks that points in equation~\eqref{embedding}
up to terms of order $O\big(\epsilon^2\big)$ equal
 \begin{gather*}
 \left[1- \frac{7}{30} \epsilon : \frac{J_{10}^{(0)}}{J_4^{(0)} J_6^{(0)}} - \frac{2}{2025} \epsilon : 0 : 0 : -\frac{512}{9} \epsilon: 0 \right] .
\end{gather*}
This means that under the usual morphism to $\mathbb{P}^5$ the regular genus-two curves with $J_2=0$ and constant ratio $J_{10}/(J_4J_6)$ are mapped to the same point.
\end{rem}

\subsection{Recovering the equation of the curve from invariants}
Let $\p \in \M_2$ and $\mathcal{C}$ a genus-two curve corresponding to $\p$ def\/ined by the sextic polynomial $f$ in equation~\eqref{genus_two_curve}.
Then, $\Aut (\p)$ is a f\/inite group as described in \cite{SV}. The quotient space $\mathcal{C}/\Aut (\p)$ is a genus zero curve and therefore isomorphic
to a conic. Since conics are in one to one correspondence with three-by-three symmetric matrices (up to equivalence), let $M=\left[ A_{i j} \right]$
be the symmetric matrix corresponding to this conic. Let $\mathbf{X} = [X_1: X_2 : X_3] \in \P^2$ and
\begin{gather}\label{conic-1}
 \mathcal Q \colon \  \mathbf{X}^t \cdot M \cdot \mathbf{X} = \sum_{i,j=1}^3 A_{ij} X_i X_j = 0.
\end{gather}
 Clebsch \cite{MR0485106} determined the entries of this matrix $M$ as follows
\begin{gather}
 A_{11} = 2 C + \frac{1}{3} A B , \nonumber\\
 A_{22} = A_{13} = D ,\nonumber\\
 A_{33} = \frac{1}{2} B D + \frac{2}{9} C \big(B^2+ A C\big) ,\nonumber\\
 A_{23} = \frac{1}{3} B \big(B^2+ A C\big) + \frac{1}{3} C \left(2 C + \frac{1}{3} A B\right) ,\nonumber\\
 A_{12} = \frac{2}{3} \big(B^2+ A C\big) .\label{Clebsch_invariants2}
\end{gather}

The coef\/f\/icients are obtained as follows: from the sextic $f$ in equation~\eqref{genus_two_curve} three binary quadrics
$\mathsf{y}_i(x) $ with $i=1,2,3$ are obtained by an operation called `\"Uberschiebung' \cite[p.~317]{Me} or transvection. The quadrics $\mathsf{y}_i$ for $i=1,2,3$ have the property that their coef\/f\/icients are polynomial expressions in the coef\/f\/icients of $f$ with rational coef\/f\/icients. Moreover, under the operation $f(x) \mapsto \tilde{f}(x)=f(-x)$ the quadrics change according to $\mathsf{y}_i(x) \mapsto \tilde{\mathsf{y}}_i(x)= \mathsf{y}_i(-x)$ for $i=1,2,3$. Hence, they are not invariants of the sextic $f$. The coef\/f\/icients $A_{ij}$ in equation~\eqref{Clebsch_invariants2} satisfy $A_{ij}=(\mathsf{y}_i\mathsf{y}_j)_2$.\footnote{For two binary forms $f$, $g$ of degree $m$ and $n$, respectively, we denote the \"Uberschiebung of order $k$ by $(fg)_k=(-1)^k (gf)_k$. For $\tilde{f}(x)=f(-x)$ and $\tilde{g}(x)=g(-x)$ and $m=n=k$, we have $(fg)_m
=(-1)^m (\tilde{f} \tilde{g})_m$.} Therefore, the coef\/f\/icients $A_{ij}$ are invariant under the operation $f(x) \mapsto \tilde{f}(x)=f(-x)$, and
the locus $D=0$ is equivalent to
\begin{gather*}
 D = 0 \quad \Leftrightarrow \quad (\mathsf{y}_1\mathsf{y}_3)_2 = (\mathsf{y}_2\mathsf{y}_2)_2 =0 .
\end{gather*}

We def\/ine $R$ to be $1/2$ times the determinant of the three binary quadrics $ \mathsf{y}_i$ for $i=1,2,3$
with respect to the basis $x^2$, $x$, $1$. If one extends the operation of \"Uberschiebung by product rule \cite[p.~317]{Me},
then $R$ can be re-written as
\begin{gather*}
 R = - (\mathsf{y}_1\mathsf{y}_2)_1 (\mathsf{y}_2\mathsf{y}_3)_1 (\mathsf{y}_3\mathsf{y}_1)_1 ,
\end{gather*}
or, equivalently, as
\begin{gather*}
 R = -\frac{1}{8} \big( \mathsf{y}_{1,yy} \mathsf{y}_{2,xy} \mathsf{y}_{3,xx} - \mathsf{y}_{1,yy} \mathsf{y}_{2,xx} \mathsf{y}_{3,xy} - \mathsf{y}_{1,xy} \mathsf{y}_{2,yy} \mathsf{y}_{3,xx}\\
 \hphantom{R=}{}+  \mathsf{y}_{1,xy} \mathsf{y}_{2,xx} \mathsf{y}_{3,yy} + \mathsf{y}_{1,xx} \mathsf{y}_{2,yy} \mathsf{y}_{3,xy} - \mathsf{y}_{1,xx} \mathsf{y}_{2,xy} \mathsf{y}_{3,yy} \big) .
\end{gather*}
It is then obvious that under the operation $f(x) \mapsto \tilde{f}(x)=f(-x)$
the determinant $R$ changes its sign, i.e., $R(f) \mapsto R(\tilde{f})=-R(f)$.
A straightforward calculation shows that
\begin{gather*}
 R^2 = \frac{1}{2} \left| \begin{matrix} A_{11} & A_{12} & A_{13} \\ A_{12} & A_{22} & A_{23} \\ A_{13} & A_{23} & A_{33} \end{matrix}\right| ,
\end{gather*}
where $A_{ij}$ are the invariants in equation~\eqref{Clebsch_invariants2}.
Like the coef\/f\/icients $A_{ij}$, $R^2$ is invariant under the operation
$f(x) \mapsto \tilde{f}(x)=f(-x)$ and must be a polynomial in $(J_2, J_4, J_6, J_{10})$.
Substitu\-ting~\eqref{invariants} into the Clebsch invariants and then equation~\eqref{Clebsch_invariants2} it follows that
\begin{gather}\label{Rsqr}
R^2 = \left( 2^{9} 3^{-9} 5^{-10} i \dfrac{\chi_{35}(\underline{\tau})}{\chi_{10}(\underline{\tau})^2} \right)^2 .
\end{gather}

Bolza \cite{MR1505464} described the possible automorphism groups of genus-two curves def\/ined by sextics and provided criteria for
the cases when the automorphism group of the sextic curve in equation~\eqref{genus_two_curve} is nontrivial. For a detailed discussion of the automorphism groups of genus-two curve def\/ined over any f\/ield $k$ and the corresponding loci in $\M_2$ see \cite{SV}. We have the following lemma summarizing our discussion:

\begin{lem} We have the following statements:
\begin{enumerate}\itemsep=0pt
\item[$1.$] $R^2$ is an order $30$ invariant of binary sextics expressed as a polynomial in $(J_2, J_4, J_6, J_{10})$ as in {\rm \cite[equation~(17)]{SV}} given by plugging Clebsch invariants and~\eqref{Clebsch_invariants2} into equation~\eqref{Rsqr}.
\item[$2.$] The locus of curves $\p \in \M_2$ such that $V_4 \embd \Aut (\p)$ is a two-dimensional irreducible rational subvariety of $\M_2$ given by the equation $R^2 = 0$ and a birational parametrization given by the $u,v$-invariants as in {\rm \cite[Theorem~1]{SV}}.
\end{enumerate}
\end{lem}

We have introduced the invariant $R^2$ for any binary sextic $f$. To the corresponding symmetric matrix $M$ with coef\/f\/icients $A_{ij}=(\mathsf{y}_i\mathsf{y}_j)_2$ of order zero and invariant under the operation $f(x) \mapsto \tilde{f}(x)=f(-x)$, we associated a conic $\mathcal{Q}$. Similarly, there is also a cubic curve given by the equation
\begin{gather}\label{cubic}
\mathcal T\colon \ \sum_{1 \leq i, j, k \leq 3} a_{ijk} X_i X_j X_k=0 ,
\end{gather}
where the coef\/f\/icients $a_{ijk}$ are of order zero and invariant under $f(x) \mapsto \tilde{f}(x)=f(-x)$.
In terms of `\"Uberschiebung' the coef\/f\/icients are obtained by
\begin{gather*}
 a_{ijk} = (f\mathsf{y}_i)_2 (f\mathsf{y}_j)_2 (f\mathsf{y}_k)_2 .
\end{gather*}
The coef\/f\/icients $a_{ijk}$ are given explicitly as follows:
\begin{gather*}
 36 a_{111} = 8\big(A^2C-6BC+9D\big),\\
 36 a_{112} = 4\big(2B^3+4ABC+12C^2+3AD\big),\\
 36 a_{113} = 36 a_{122} = 4\big(AB^3 + 4/3 A^2 BC + 4 B^2 C + 6 AC^2 + 3 BD\big),\\
 36 a_{123} = 2\big(2B^4+4AB^2C+4/3 A^2 C^2 + 4 BC^2+ 3 ABD+ 12 CD\big),\\
 36 a_{133} = 2\big(AB^4 + 4/3 A^2B^2C+ 16/3 B^3 C  + 26/3 ABC^2+ 8 C^3+ 3 B^2 D + 2 ACD\big),\\
 36 a_{222} = 4\big( 3B^4 + 6 AB^2 C + 8/3 A^2 C^2 + 2 BC^2- 3CD\big),\\
 36 a_{223} = 2\big({}-2/3 B^3 C- 4/3 ABC^2- 4C^3 + 9 B^2D+8ACD\big) ,\\
 36 a_{233} = 2\big(B^5+2AB^3C+8/9 A^2BC^2+ 2/3 B^2 C^2 - BCD + 9 D^2\big),\\
 36 a_{333} = -2 B^4C - 4AB^2C^2- 16/9 A^2 C^3 - 4/3 BC^3  + 9 B^3 D + 12 ABCD + 20 C^2 D .
\end{gather*}
The relations between all aforementioned invariants and Siegel modular forms, in particular the relation between $\chi_{35}$ and $R^2$ can be found in \cite{data}.

Since `\"Uberschiebung' preserves the rationality of the coef\/f\/icients, we have the following corollary:
\begin{cor}
Let $\p \in \M_2$ and $\mathcal{C}$ a genus-two curve corresponding to $\p$ defined by a sextic polynomial $f$ in equation~\eqref{genus_two_curve}.
Then, $\Aut (\p)$ is a finite group, and the quotient space $\mathcal{C}/\Aut (\p)$ is a genus zero curve isomorphic to the conic $\mathcal Q$ in equation~\eqref{conic-1}. Moreover, if $\p \in \M_2 (K)$, for some number field $K$, the conic $\mathcal Q$ and cubic $\mathcal T$ have $K$-rational coefficients.
\end{cor}

The intersection of the conic $\mathcal Q$ with the cubic $\mathcal T$ consists of six points which are the zeroes of a polynomial $f(x)$  of degree 6 in the parameter $x$. The roots of this polynomial are the images of the Weierstrass points under the hyperelliptic projection. Hence, the af\/f\/ine equation of a~genus-two curve corresponding to $\p$ is given by $y^2= f(x)$. The main question is if the sextic given by $y^2= f(x)$ provides a genus-two curve def\/ined over a minimal f\/ield of def\/inition. We start with the following known result.

\begin{prop}\label{Mestre}
A genus $g\geq 2$ hyperelliptic curve $\X_g$ with hyperelliptic involution $w$ is defined over the $K$ if and only if the conic $\mathcal Q = \X_g/\< w\>$ has a $K$-rational point.
\end{prop}

The above result was brief\/ly described in \cite[Lemma 1]{Me} even though it seems as it had been known before.
Mestre's method is brief\/ly described as follows: if the conic $\mathcal Q$ has a rational point over $\Q$, then this leads to a parametrization of $\mathcal Q$, say $(h_1 (x), h_2( x), h_3(x) )$. Substitute $X_1$, $X_2$, $X_3$ by $h_1 (x)$, $h_2( x)$, $h_3(x)$ in the cubic $\mathcal T$ and we get the degree 6 polynomial~$f(x)$. However, if the conic has no rational point or $R^2= \frac{1}{2} \det{M} =0$ the method obviously fails. In Section~\ref{universal-curve} we determine the intersection $\mathcal T \cap \mathcal Q$ over a quadratic extension which is always possible.

\section{A universal genus-two curve from the moduli space}\label{universal-curve}

The goal of this section is to explicitly determine a universal equation of a genus-two curve corresponding to this generic point $\p$.
We have the following lemma:
\begin{lem}\label{lem-Conic-Diag}The conic $\mathcal{Q}$ in equation~\eqref{conic-1} for $J_4 \cdot J_6 \cdot J_{10} \not = 0$
is equivalent over $\mathbb{Q}[J_2,\rho,\sigma,\kappa]$ to the conic
\begin{gather}\label{ConicRenorm}
\mathcal{Q}'\colon \ x_1^2 -\gamma x_2^2 - \Lambda_6 x_3^2 =0 ,
\end{gather}
where $(\rho,\sigma,\kappa)$ are the absolute invariants in equation~\eqref{OurCoordinates}, $\gamma=\rho^2 + \sigma$ and
\begin{gather}
\Lambda_6 = -\gamma^3-27 \rho \gamma^2-81 \rho^2 (\rho+12) \gamma+729 \rho^2 (\rho+12)^2\nonumber\\
\hphantom{\Lambda_6 =}{} + \big({-}6 \rho \gamma^2+54 \rho (5 \rho+36) \gamma -1944 \rho^2 (\rho+12)\big) \kappa\nonumber\\
\hphantom{\Lambda_6 =}{}+\big(9 \gamma^2-9 \rho (\rho+36) \gamma+162 \rho (\rho^2+32 \rho+144)\big) \kappa^2\nonumber\\
\hphantom{\Lambda_6 =}{}+ ((30 \rho+216) \gamma-432 \rho (\rho+12) ) \kappa^3+\big(9 \rho^2-24 \gamma+504 \rho+1296\big) \kappa^4\nonumber\\
\hphantom{\Lambda_6 =}{}+ (-24 \rho-288 ) \kappa^5+16 \kappa^6.\label{ConicRenormCoeffs}
\end{gather}
Moreover, for $J_2,\rho,\sigma,\kappa \in \mathbb{Q}$ the conic $\mathcal{Q}$ in equation~\eqref{conic-1} has a rational point
if and only if the conic $\mathcal{Q}'$ in equation~\eqref{ConicRenorm} does.
\end{lem}

\begin{proof} For the conic $\mathcal{Q}$ in equation~\eqref{conic-1}, we apply the coordinate transformation given by
\begin{gather}
 X_1 = 2  ( AB+6C  )^{4} \big( AC+B^{2} \big) x_1\nonumber\\
\hphantom{X_1 =}{} + 108 B ( AB+6 C )^{2} \big( 4A^2 C^2+8A B^2 C+4 B^4-3ABD -18 C D \big) x_2\nonumber\\
\hphantom{X_1 =}{} + 41990 B^3 \big( 8 A^2 B C^2+14 A B^3 C+6 B^5+12 A C^3+12 B^2 C^2-27D^2 \big) x_3,\nonumber\\
 X_2 = -  ( AB+6C  )^5 x_1 -419904 B^3 x_3\nonumber\\
\hphantom{X_2 =}{} \times \big(4 A^2 B^2 C+3 A B^4+30 A B C^2+18 B^3 C-18 A C D-18 B^2 D+36 C^3\big),\nonumber\\
 X_3 = -2^6 3^9 B^3 \big(4 A^2 C^2+8 A B^2 C+4 B^4-3 A B D-18 C D\big) x_3.\label{TransfoDiag}
\end{gather}
We then obtain the conic $\mathcal{Q}'$ in equation~\eqref{ConicRenorm}.
Equation~\eqref{TransfoDiag} can be rewritten as transformation over $\mathbb{Q}[J_2,\rho,\sigma,\kappa]$
using equations~\eqref{iClebsch_invariants} and~\eqref{OurCoordinates}.
\end{proof}

We have the following lemma:

\begin{lem}\label{lem-Conic-Diag-BasePt}
Assume $\rho,\sigma,\kappa \in \mathbb{Q}$. The conic $\mathcal{Q}'$ in equation~\eqref{ConicRenorm} has a rational point
if and only if there are rational numbers $\alpha, \beta \in \mathbb{Q}$ such that
\begin{gather}\label{ConicRelation}
 \alpha^2 + \Lambda_6 \beta^2\sigma = \gamma .
\end{gather}
The rational point on the conic $\mathcal{Q}'$ is then given by
\begin{gather}\label{Basepoint}
\big[x_1^0:x_2^0:x_3^0\big] = [\alpha \rho +\gamma: \alpha + \rho: \beta\sigma ] .
 \end{gather}
 Conversely, every rational point on the conic $\mathcal{Q}'$ can be written in the form of equation~\eqref{Basepoint} for some rational numbers $\alpha, \beta \in \mathbb{Q}$ satisfying  equation~\eqref{ConicRelation}.
\end{lem}
\begin{proof}If rational numbers $\alpha$, $\beta$ exist such that equation~\eqref{ConicRelation} is satisf\/ied, then the point in equation~\eqref{Basepoint} is rational and is easily checked to be on the conic. If there is a rational point on the conic then we can choose $\beta\in \mathbb{Q}$ in equation~\eqref{Basepoint}, thus $\alpha \in \mathbb{Q}$.
\end{proof}

We have the following:
\begin{lem}\label{lem-ParamDiagY}
Assume that a point on the conic in equation~\eqref{ConicRenorm} is given by equation~\eqref{Basepoint} with
$x_2^0 \not = 0$ which is always possible if $\rho \not =0$. Then every point on the conic is given by
\begin{gather}
 x_1 = (\alpha \rho+\gamma) U^2 + 2 \Lambda_6 \beta \sigma U V+ \Lambda_6 (\alpha \rho+\gamma) V^2, \nonumber\\
 x_2 = (\alpha+\rho) U^2- \Lambda_6 (\alpha+\rho) V^2 ,\nonumber\\
 x_3 = \beta\sigma U^2+2 (\alpha \rho+\gamma) U V +\Lambda_6 \beta \sigma V^2 ,\label{ParamDiagY}
\end{gather}
for some $[U:V]\in \mathbb{P}^1$. The parametrization in equation~\eqref{ParamDiagY} is a rational parametrization of the conic $\mathcal{Q}'$ if and only if
$\alpha, \beta, \rho,\kappa, \sigma \in \mathbb{Q}$.
\end{lem}
\begin{proof}
If a point of $\mathcal{Q}'$ is obtained from some (rational) values $(\alpha,\beta)$ then there are three more (rational)
points given by setting $(\alpha,\beta) \mapsto (\pm \alpha, \pm \beta)$.
If $\rho \not=0$, one of these points satisf\/ies $x_2^0 =\alpha + \rho \not = 0$.
The proof then follows from the known formulas parametrizing conics for $x_2^0 \not =0$ given by
\begin{gather*}
x_1= a x_1^0 U^2-2 c x_3^0 U V - c x_1^0 V^2, \\
x_2= a x_2^0 U^2+c x_2^0 V^2, \\
x_3= a x_3^0 U^2+2 a x_1^0 U V-c x_3^0 V^2 ,
\end{gather*}
where $a=1$, $b=-\gamma$, $c=-\Lambda_6$ and $x_1^0$, $x_2^0$, $x_3^0$ were given in equation~\eqref{Basepoint}.
\end{proof}

\begin{rem}
If $\alpha = \rho = 0$ and $\gamma \not = 0$, a formula similar to equation~\eqref{ParamDiagY}
can be found using the fact that $x_1^0 \not = 0$ in equation~\eqref{Basepoint} in this case.
\end{rem}

\begin{rem}
If a point of $\mathcal{Q}'$ is obtained for some (rational) values $(\alpha,\beta)$ then three more (rational)
points on $\mathcal{Q}'$ are given by setting $(\alpha,\beta) \mapsto (\pm \alpha, \pm \beta)$ in equation~\eqref{Basepoint}.
\end{rem}

Changing from coordinates $[X_1:X_2:X_3]$ to coordinates $[x_1:x_2:x_3]$ transforms the conic~$\mathcal{Q}$ in equation~\eqref{conic-1}
into the conic $\mathcal{Q}'$ in equation~\eqref{ConicRenorm}. Similarly, under the same change of coordinates the cubic $\mathcal{T}$
in equation~\eqref{cubic} becomes
\begin{gather}
\Lambda_1 (18\gamma+\Lambda_3) x_1^3+ \gamma^3 \Lambda_2 x_2^3 - (\gamma - \Lambda_1) \Lambda_6^2 x_3^3+ 3 \Lambda_1 \Lambda_6 x_1^2 x_3 + 3 \gamma (9\gamma\rho+ \kappa \Lambda_3) x_1^2 x_2 \nonumber\\
\hphantom{\mathcal{T}'\colon \  0=}{} + 3 \gamma^2 \Lambda_3 x_1 x_2^2 + 3 \Lambda_5 \Lambda_6 x_1 x_3^2 + 3 \gamma \Lambda_4 \Lambda_6 x_2 x_3^2+ 3 \gamma^2 \Lambda_6 x_2^2 x_3
+ 6 \gamma \kappa \Lambda_6 x_1x_2x_3,\label{CubicRenorm}
\end{gather}
with coef\/f\/icients given by
\begin{gather}
\Lambda_1 = 9\rho+\kappa^2,\nonumber\\
\Lambda_2  = \gamma+18\rho+3\rho\kappa-4\kappa^2,\nonumber\\
\Lambda_3  = 27\rho(\rho+12)+(\gamma-36\rho)\kappa+3(\rho+12)\kappa^2-4\kappa^3,\nonumber\\
\Lambda_4 = -\gamma^2-9\gamma\rho-3(\gamma\rho-9\rho(\rho+12))\kappa+(5\gamma-36\rho)\kappa^2+3(\rho+12)\kappa^3-4\kappa^4,\nonumber\\
\Lambda_5 = -27\gamma\rho(\rho+6)+243\rho^2(\rho+12)-\big(\gamma^2-45\gamma\rho+324\rho^2\big)\kappa\nonumber\\
\hphantom{\Lambda_5 =}{} - (3\gamma(\rho-6)-54\rho(\rho+12) )\kappa^2+(5\gamma-72\rho)\kappa^3+3(\rho+12)\kappa^4-4\kappa^5.\label{CubicRenormCoeffs}
\end{gather}

We also discuss the conic, cubic, rational point and parametrization in the cases where $J_2=0$ and $J_4 \cdot J_6=0$:
\begin{lem}\label{lem-special-cases}
If $J_2 =J_4=0$ and $J_6\cdot J_{10}\neq 0$, the conic $\mathcal{Q}$ in equation~\eqref{conic-1} is equivalent over $\mathbb{Q}[J_6, J_{10}]$ to the conic
\begin{gather*}
\mathcal{Q}'\colon \ x_1^2 -\mu x_2^2 -(1-\mu) x_3^2 =0 ,
\end{gather*}
with $\mu=J_6^5/\big(2^4 3^4 5^5 J_{10}^3\big)$ and a rational point given by $\big[x_1^0:x_2^0:x_3^0\big] =[1:1:1]$. A rational parametrization of $\mathcal{Q}'$ is then given by
\begin{gather*}
x_1 = U^2+2(1-\mu)UV +(1-\mu)V^2, \qquad x_2 = U^2-(1-\mu)V^2, \\
 x_3 = U^2+2UV+(1-\mu)V^2
\end{gather*}
with $[U:V]\in \mathbb{P}^1$. Under the same change of coordinates the cubic $\mathcal{T}$ in equation~\eqref{cubic} becomes
\begin{gather*}
\mathcal{T}'\colon  \  0 = 2 x_1^3 - \mu^2 x_2^3-2 (1-\mu)^2 x_3^3 - 6 \mu x_1^2 x_2 \\
\hphantom{\mathcal{T}'\colon  \  0 =}{} - 6 (1-\mu) x_1^2 x_3
+ 6(1-\mu) x_1 x_3^2 - 3\mu(1-\mu) x_2 x_3^2.
\end{gather*}

If $J_2 =J_6=0$ and $J_4 \cdot J_{10}\neq 0$, the conic $\mathcal{Q}$ in equation~\eqref{conic-1} is equivalent over
$\mathbb{Q}[J_4, J_{10}]$ to the conic
\begin{gather*}
\mathcal{Q}'\colon \ x_1^2 - x_2^2 - (1-\nu) x_3^2 =0 ,
\end{gather*}
with $\nu=J_4^5/\big(2^2 3^5 5^5 J_{10}^2\big)$ and a rational point given by $\big[x_1^0:x_2^0:x_3^0\big] =[1:1:0]$. A rational parametrization of $\mathcal{Q}'$ is then given by
\begin{gather*} 
x_1 = U^2+(1-\nu)V^2, \qquad x_2 = U^2-(1-\nu)V^2, \qquad x_3 = 2 UV
\end{gather*}
with $[U:V]\in \mathbb{P}^1$. Under the same change of coordinates the cubic $\mathcal{T}$ in equation~\eqref{cubic} becomes
\begin{gather*}
\mathcal{T}'\colon \  0 = \big(1-\nu^2\big) x_1^3 - \nu^2 x_2^3-(1-\nu)^2 x_3^3 + \nu(1-3\nu) x_1^2 x_2 \\
\hphantom{\mathcal{T}'\colon \  0 =}{} - (1-\nu) (3+\nu) x_1^2 x_3 + \nu(1-3\nu) x_1 x_2^2 -\nu (1-\nu) x_2^2 x_3 \\
\hphantom{\mathcal{T}'\colon \  0 =}{}+ (1-\nu)(3-\nu) x_1 x_3^2 + \nu(1-\nu)x_2x_3^2 -2\nu(1-\nu)x_1x_2x_3.
\end{gather*}
\end{lem}
\begin{proof}The proof is analogous to the proofs of Lemmas~\ref{lem-Conic-Diag-BasePt} and \ref{lem-ParamDiagY}.
\end{proof}

\begin{rem}The absolute invariants $(\rho, \sigma, \kappa)$ in equation~\eqref{OurCoordinates} such that $\rho$ and $\kappa$ do not vanish simultaneously
and $J_{10} \neq 0$ describe the moduli of genus-two curves with $J_4 \cdot J_6 \cdot J_{10} \neq 0$. The discussion of Lemma~\ref{lem-special-cases} proves that only for genus-two curves with $J_4 \cdot J_6 \cdot J_{10} \neq 0$, the conic~$\mathcal{Q}$ in equation~\eqref{conic-1} is \emph{not} guaranteed to have a rational point.
\end{rem}

Substituting the parametrization of the conic $\mathcal{Q}'$ in Lemma~\ref{lem-ParamDiagY} into the cubic $\mathcal{T}'$ in equation~\eqref{CubicRenorm} and setting $U=x$ and $V=1$, one obtains the ramif\/ication locus of a sextic curve. The ramif\/ication locus is equivalent to $f(x) =\sum\limits_{i=0}^6 a_i(\alpha,\beta) x^i=0$, where we write the sextic polynomial in the form
\begin{gather}
f(x) = \big(d^{(1)}_0 + d^{(2)}_0\big) x^6
 + \big(d^{(1)}_1 + d^{(2)}_1\big) \Lambda_6x^5
 + \big(d^{(1)}_2 + d^{(2)}_2\big) \Lambda_6x^4 \nonumber\\
\hphantom{f(x)=}{} + d^{(1)}_3 \Lambda^2_6x^3
 + \big(d^{(1)}_2 - d^{(2)}_2\big)\Lambda_6^2 x^2
 + \big(d^{(1)}_1 - d^{(2)}_1\big) \Lambda_6^3 x
 + \big(d^{(1)}_0 - d^{(2)}_0\big) \Lambda^3_6.\label{SexticPolynomial}
\end{gather}
In terms of the coordinates of the point $\big[x_1^0: x_2^0: x_3^0\big]$ in equation~\eqref{Basepoint} we have set
\begin{gather}
 d^{(1)}_{j} = c^{(1)}_{j,0} \big[\big(\rho^2+\gamma\big) x_1^0 + 2\gamma\rho x_2^0\big]
 			\gamma		+ 3 c^{(1)}_{j,1} \big[x_1^0+ \rho x_2^0\big]
			\gamma \Lambda_6^{\delta_j} x_3^0 \nonumber\\
\hphantom{d^{(1)}_{j} =}{}- c^{(1)}_{j,0} \big[\big(\rho^2+\gamma\big)x_1^0 - 2\gamma\rho x_2^0\big] \rho^2\sigma^{-2}\Lambda_6 \big(x_3^0\big)^2
 	- 3 c^{(1)}_{j,1} \big[x_1^0- \rho x_2^0 \big] \rho^2 \sigma^{-2} \Lambda_6^{1+\delta_j}\big(x_3^0\big)^3 \nonumber\\
\hphantom{d^{(1)}_{j} =}{}+ \big[3 \sigma c^{(1)}_{j,2} -2 \rho^2c^{(1)}_{j,0}\big] \sigma^{-1} \Lambda_6 x_1^0 \big(x_3^0\big)^2 + c^{(1)}_{j,3}\Lambda_6^{1+\delta_j} \big(x_3^0\big)^3 ,\nonumber\\
 d^{(2)}_{j} = c^{(2)}_{j,0} \big[2\rho x_1^0+ \big(\rho^2+\gamma\big)x_2^0 \big]\gamma^2
 	+ 6 c^{(2)}_{j,1} \big[\rho x_1^0 + \gamma x_2^0 \big] \gamma \Lambda_6^{\delta_j} x_3^0 \nonumber\\
\hphantom{d^{(2)}_{j} =}{}- c^{(2)}_{j,0} \big[ 2 \gamma x_1^0 - \rho\big(\rho^2+\gamma\big)x_2^0\big] \gamma \rho \sigma^{-2} \Lambda_6 \big(x_3^0\big)^2
 	-6 c^{(2)}_{j,1} \big[x_1^0 - \rho x_2^0\big] \gamma\rho\sigma^{-2} \Lambda_6^{1+\delta_j} \big(x_3^0\big)^3 \nonumber\\
 \hphantom{d^{(2)}_{j} =}{}+ c^{(2)}_{j,2} \gamma \Lambda_6^{\delta_j} x_1^0 \big(x_3^0\big)^2.\label{SexticPolynomialCoeffsDecomp}
\end{gather}
All coef\/f\/icients remain regular and in general non-vanishing for $\sigma=0$ since
$x_0^3/\sigma=\beta$. Here, $(\alpha,\beta)$ is a pair solving equation~\eqref{ConicRelation}, and
the coef\/f\/icients $c^{(n)}_{j,k}$ are given by
\begin{alignat*}{3}
& c^{(1)}_{0,0}	= 18\gamma\Lambda_1+3 \gamma \Lambda_3+\Lambda_1\Lambda_3, \qquad && c^{(1)}_{0,1} 	= \gamma+\Lambda_1,&\\
& c^{(1)}_{0,2} = -\gamma \Lambda_3+\Lambda_5,\qquad && c^{(1)}_{0,3} = -4\gamma+\Lambda_1,	&\\
& c^{(2)}_{0,0} = 27 \gamma\rho+\gamma\Lambda_2+3 \kappa\Lambda_3, \qquad && c^{(2)}_{0,1}	= \kappa,&\\
& c^{(2)}_{0,2}	= -\gamma\Lambda_2+3\Lambda_4,\qquad &&&\\
& c^{(1)}_{1,0}	= 6(\gamma+\Lambda_1), \qquad && c^{(1)}_{1,1}	= 2(18\gamma\Lambda_1+\gamma\Lambda_3+\Lambda_1\Lambda_3+2\Lambda_5),&\\
 &c^{(1)}_{1,2}	=2(-2\gamma+3\Lambda_1),\qquad && c^{(1)}_{1,3}	= 6(-\gamma\Lambda_3+\Lambda_5),&\\
& c^{(2)}_{1,0}	=12\kappa,\qquad && c^{(2)}_{1,1}	=18\gamma\rho+2\kappa\Lambda_3+2\Lambda_4,&\\
& c^{(2)}_{1,2}	= 12\kappa\Lambda_6,\qquad &&&\\
& c^{(1)}_{2,0}	= 3(18\gamma\Lambda_1-\gamma\Lambda_3+\Lambda_1\Lambda_3+4\Lambda_5), \qquad && c^{(1)}_{2,1}	= 5(-\gamma+3\Lambda_1),&\\
& c^{(1)}_{2,2}	= 72\gamma\Lambda_1+\gamma\Lambda_3+4\Lambda_1\Lambda_3+11\Lambda_5, \qquad && c^{(1)}_{2,3}	=15 \Lambda_1,&\\
& c^{(2)}_{2,0}	=3(9\gamma\rho-\gamma\Lambda_2+\kappa\Lambda_3+4\Lambda_4), \qquad && c^{(2)}_{2,1}	=5\kappa,&\\
& c^{(2)}_{2,2}	=3(36\gamma\rho+\gamma\Lambda_2+4\kappa\Lambda_3+\Lambda_4),\qquad && &\\
& c^{(1)}_{3,0}	= 20(-\gamma+\Lambda_1),\qquad && c^{(1)}_{3,1}	=4(18\gamma\Lambda_1- \gamma \Lambda_3+\Lambda_1\Lambda_3+4\Lambda_5),&\\
& c^{(1)}_{3,2}	=20 \Lambda_1,\qquad && c^{(1)}_{3,3}	=4(36\gamma\Lambda_1+3 \gamma \Lambda_3+2\Lambda_1\Lambda_3+3\Lambda_5).&
\end{alignat*}
The coef\/f\/icients $\Lambda_1, \dots, \Lambda_5$ and $\Lambda_6$ were given in equations~\eqref{CubicRenormCoeffs}
and~\eqref{ConicRenormCoeffs}, respectively.

\begin{rem}Equation~\eqref{SexticPolynomialCoeffsDecomp} allows to easily describe the change in the sextic polynomial under the action of the automorphism of the conic~$\mathcal{Q}'$ given by $[x_1:x_2:x_3] \mapsto [\pm x_1: \pm x_2: x_3]$.
\end{rem}

We make the following remarks:

\begin{rem}The transformation $x\to \frac {\Lambda_6} {x}$ maps the coef\/f\/icients $d^{(1)}_k \pm d^{(2)}_k \mapsto d^{(1)}_k \mp d^{(2)}_k$ for $k=0,1,2$ and $a_3 \mapsto a_3$. This is to be expected since the coef\/f\/icients are in terms of invariants of the binary sextic $f(x, z)$ and $x \to \frac 1 x$ just permutes $x$ and $z$.
\end{rem}

\begin{rem}The fact that the coef\/f\/icients def\/ining $\mathcal{C}_{( \alpha, \beta)}$ are polynomials of the new absolute invariants $\gamma$, $\kappa$, $\sigma$ and appear in the particular pattern given by equation~(\ref{SexticPolynomialCoeffsDecomp}) opens up the question about their meaning. It turns out that in the context of the F-theory/heterotic string theory duality the new invariants parameterize the physical defects of a certain class of 6d $\mathcal{N}=(1,0)$ non-geometric vacua of the heterotic string when dualizing to F-theory. We will address this question in more detail in future work.
\end{rem}

We have the following main result:
\begin{thm}\label{main-thm}Let $\p \in \M_2$ such that $\p \in \M_2 (K)$, for some number field $K$, and $\mathfrak j=[J_2 : J_4 : J_6 : J_{10} ]$ the corresponding point in $\wM (\mathcal O_K)$, where $\mathcal O_K$ is the ring of integers of~$K$. A~genus-two curve corresponding to~$\p$ is constructed as follows:
\begin{enumerate}\itemsep=0pt
\item[$i)$] If $J_2 \cdot J_{10}\neq 0$ there is a genus-two curve $\mathcal{C}_{(\alpha,\beta)}$
given by
\begin{gather}\label{Sextic_main}
 \mathcal{C}_{(\alpha, \beta)}\colon \ y^2= \sum_{i=0}^6 a_{i}(\alpha,\beta) x^i ,
 \end{gather}
with coefficients given in equations~\eqref{SexticPolynomial} and~\eqref{SexticPolynomialCoeffsDecomp}, and a pair $(\alpha,\beta)$ satisfying
\begin{gather*} \alpha^2 + \Lambda_6 \beta^2\sigma = \gamma, \end{gather*}
where $\Lambda_6$, $\sigma$, and $\gamma$ are determined by $\p$. Moreover, $\mathcal{C}_{(\alpha, \beta)}$ is defined over its field of moduli~$K$, i.e., $a_{i}(\alpha,\beta) \in K$, $i=0, \dots , 6$, if and only if $K$-rational~$\a$ and~$\beta$ exist.

\item[$ii)$] If $J_2=0$ and $J_4\cdot J_6 \cdot J_{10}\neq 0$, there is a genus-two curve given by setting $\rho=\kappa\neq0$ in equation~\eqref{Sextic_main}.

\item[$iii)$] If $J_2 =J_6=0$ and $J_4 \cdot J_{10}\neq 0$, there is only one genus-two curve given by
\begin{gather*}
y^2 =  (4\nu+1)(2\nu-1) x^6+2(1-\nu)(4\nu+3)x^5-15(1-\nu)x^4\\
\hphantom{y^2 =}{} +20(1-\nu)^2 x^3+5 (2\nu-3) (1-\nu)^2 x^2+6 (1-\nu)^3 x -(1-\nu)^3
\end{gather*}
with $\nu=J_4^5/\big(2^2 3^5 5^5 J_{10}^2\big)$.

\item[$iv)$] If $J_2 =J_4= 0$ and $J_6 \cdot J_{10}\neq 0$, there is only one genus-two curve given by
\begin{gather*}
y^2 =  5 x^6+12(1-\mu)x^5-15(1-\mu)x^4-80(1-\mu)^2 x^3\\
\hphantom{y^2 =}{} +15 (4\mu-7) (1-\mu)^2 x^2 -60 (1-\mu)^3 x+(4\mu-13)(1-\mu)^3.
\end{gather*}
with $\mu=J_6^5/\big(2^4 3^4 5^5 J_{10}^3\big)$.

\item[$v)$] If $J_2 =J_4=J_6= 0$ and $J_{10} \neq 0$, there is only one genus-two curve given by
\begin{gather*} y^2 = x^6-x .\end{gather*}
\end{enumerate}
\end{thm}

\begin{proof} We already proved that there are genus-two curves $y^2=f(x)$ corresponding to $\p$, where~$f(x)$ is given in equation~\eqref{SexticPolynomial}. We obtain coef\/f\/icients $ a_i(\alpha,\beta) \in \Q [\alpha,\beta,J_2, J_4, J_6, J_{10}]$ for $0\le i\le6$. The f\/ield of moduli $K$ of the point $\p$ is $K=\Q (\x_1, \x_2, \x_3)$.  For $J_2 \not= 0$ the invariants $(\rho,\sigma,\kappa)$ are birationally equivalent to $(\x_1,\x_2,\x_3)$ over $\mathbb{Q}$ by Lemma~\ref{lem-NewInv}. By Lemma~\ref{lem-Conic-Diag} the conic~$\mathcal{Q}$ in equation~\eqref{conic-1} had a $K$-rational point if and only if the conic $\mathcal{Q}'$ in equation~\eqref{ConicRenorm} does. By Lemma~\ref{lem-Conic-Diag-BasePt} the conic $\mathcal{Q}'$ has a~$K$-rational point, i.e., there is a~$K$-rational solution~$(\alpha,\beta)$ of equation~\eqref{ConicRelation}. Therefore, $a_i( \alpha,\beta) \in K$, for $i=0, \dots , 6$. The cases with $J_4 \cdot J_6=0$ are similarly obtained by applying Lemmas~\ref{lem-NewInv} and~\ref{lem-special-cases}.

 This completes the proof.
\end{proof}

\begin{rem} The four pairs $(\pm \alpha, \pm \beta)$ belong to the same conic $\mathcal{Q}'$. Therefore, we get four genus-two curves in Theorem~\ref{main-thm}, but they are all twists of each other. That is, we get one curve (over the algebraic closure), but four twists.
\end{rem}

The main benef\/it of the above result is that it will give a curve def\/ined over $\Q$ whenever possible. This is an improvement from results in~\cite{Me} where a curve is provided only for curves with automorphism group of order 2 and $J_2\neq 0$. The equation is valid even when the f\/ield of moduli is not a f\/ield of def\/inition. Hence, for every point $\p \in \M_2$ we get a curve.
Next we have the following result:
\begin{cor}\label{main-cor}For every point $\p \in \M_2$ such that $\p \in \M_2 (K)$, for some number field $K$, there is a genus-two curves $\mathcal{C}$ given by
\begin{gather*} \mathcal{C}_{(\alpha,0)}\colon \ y^2= \sum_{i=0}^6 a_{i}(\alpha,0) x^i , \end{gather*}
 corresponding to $\p$, such that $a_{i}(\alpha,0) \in K( \alpha )$, $i=0, \dots , 6$ as given in equation~\eqref{SexticPolynomial}.
 Moreover, $\mathcal{C}_{(\alpha,0)}$ is at worst defined over the quadratic extension $K ( \alpha )$ of the field of moduli $K$ with $\alpha^2=\rho^2+\sigma$.
\end{cor}

We have the immediate consequence:

\begin{cor}Let $\x_1$, $\x_2$, $\x_3$ be transcendentals. There exists a genus-two curve $\mathcal{C}_{(\alpha,0)}$ defined over $\Q (\x_1, \x_2, \x_3) [\alpha]$
with $\alpha^2=\rho^2+\sigma$
such that
\begin{gather*} \x_1(\mathcal{C}_{(\alpha,0)})=\x_1, \qquad \x_2 (\mathcal{C}_{(\alpha,0)}) = \x_2, \qquad \x_3 (\mathcal{C}_{(\alpha,0)}) = \x_3. \end{gather*}
\end{cor}

We have the following corollary:

\begin{cor}
Let $\sigma = 0$ and $\rho \not =0$ for $\p \in \M_2$. Then, there is a genus-two curve $\mathcal{C}$ given by Corollary~{\rm \ref{main-cor}}, and it is defined over the field of moduli.
\end{cor}
\begin{proof}
For $\sigma=0$ and $\rho \not =0$, we have $\gamma = \rho^2$, and we choose the $K$-rational solution $(\alpha,\beta)=(\rho,0)$ in  equation~\eqref{ConicRelation}.
\end{proof}

\begin{rem}It is easy to check using equation~\eqref{OurCoordinates} that the locus $\sigma = 0$ and $\rho \not =0$ for $\p \in \M_2$ corresponds to the locus
\begin{gather*}
J_{10} = - 2^{-11} 3^{-3} 5^{-5} \big(9 J_2^5-700 J_2^3 J_4+2400 J_2^2 J_6-262400 J_2 J_4^2+768000 J_4 J_6\big) .
\end{gather*}
\end{rem}

We have the following lemma:
\begin{lem}
In terms of the invariants $\rho$, $\sigma$, $\kappa$ and $\gamma=\rho^2+\sigma$, we have
\begin{gather*}
 D =-\frac{J_2^5 \big((\kappa-\rho)^2+9\rho\big) \big((2\kappa-\rho)^2-\gamma\big)}{2^{17} 3^7 5^5(\kappa-\rho)^5}, \\
 R^2=\frac{J_2^{15} \big((\kappa-\rho)^2+9\rho\big)^3 \Lambda_6}{2^{54}3^{21}5^{15}(\kappa-\rho)^{15}}.
\end{gather*}
In particular, the locus $D=0$ and $\chi_{35} \not =0$ is given by $\gamma=(2\kappa-\rho)^2$ or, equivalently, $\sigma=4\kappa(\kappa-\rho)$.
\end{lem}

We have the following corollary:
\begin{cor}
Let $D=0$ and $\chi_{35}^2\not =0$ for $\p \in \M_2$. Then, there is a genus-two curve $\mathcal{C}$ given by Corollary~{\rm \ref{main-cor}},
and it is defined over the field of moduli.
\end{cor}
\begin{proof}
For $\gamma=(2\kappa-\rho)^2$ we can choose $(\alpha,\beta)=(\rho-2\kappa,0)$ in equation~\eqref{ConicRelation}. As $\kappa-\rho\not= 0$ we have $y_0 \not = 0$ in equation~\eqref{Basepoint}.
\end{proof}

\subsection{A word about extra automorphisms}
In this section we derive a sextic polynomial for the sublocus of $\M_2$ with $\chi_{35}=0$. We have the following proposition:
\begin{prop}Let $D\not =0$ and $\chi_{35} =0$ for $\p \in \M_2$. Then, there is a genus-two curve $\mathcal{C}\colon y^2=F(x)$ with
\begin{gather}\label{SexticPolynomial2}
 F(x) = \big(d^{(1)}_0 + d^{(2)}_0\big) x^6
 + \big(d^{(1)}_2 + d^{(2)}_2\big) x^4
 + \big(d^{(1)}_2 - d^{(2)}_2\big) x^2
 + \big(d^{(1)}_0 - d^{(2)}_0\big),
\end{gather}
and with coefficients in $\mathbb{Z}[\alpha,\rho,\kappa]$ given by
\begin{gather}
 d^{(1)}_0  = 3\kappa\gamma^2- \big(\kappa^{2}+9\rho \big) ( 11\kappa-9\rho-126) \gamma
 		- \big(\kappa^2+9\rho \big)^{2} ( 4\kappa-3\rho-36 ),\nonumber \\
 d^{(2)}_0  =\big(\gamma^2+ \big( {-}\kappa^2+3\kappa\rho+45\rho \big) \gamma
 		-3\kappa \big(\kappa^2+9\rho \big) ( 4\kappa-3\rho-36 )\big) \alpha,\nonumber\\
d^{(1)}_2 = -15\kappa\gamma^2+15\big( \kappa^2+9\rho \big) ( 5\kappa-3\rho-18 ) \gamma
		-15\big(\kappa^2+9\rho)^2 ( 4\kappa-3\rho-36 ),\nonumber\\
 d^{(2)}_2 =\big({-}15\gamma^2+ \big( 75\kappa^2-45\kappa\rho-135\rho \big) \gamma-15\kappa
 \big(\kappa^2+9\rho \big) ( 4\kappa-3\rho-36 ) \big) \alpha.\label{SPCD2}
\end{gather}
Here, the absolute invariants $\alpha$, $\gamma$, $\rho$, $\kappa$ are subject to the constraints $\Lambda_6=0$ in equation~\eqref{ConicRenormCoeffs} and $\alpha^2=\gamma$.
\end{prop}
\begin{proof}For $\Lambda_6=\epsilon^2$ with $\epsilon \to 0$, we rescale the polynomial in equation~\eqref{SexticPolynomial2} according to $f(\epsilon x)/\epsilon^6$ before setting $\epsilon=0$. If we substitute $\Lambda_6=0$ into equation~\eqref{ConicRelation} we obtain $\alpha^2=\gamma$, \smash{$\beta=0$}. Therefore, we will use the absolute invariants $\alpha$, $\gamma$, $\rho$, $\kappa$ subject to the constraints $\Lambda_6=0$ in equation~\eqref{ConicRenormCoeffs} and $\alpha^2=\gamma$. The sextic polynomial in equation~\eqref{SexticPolynomial2} has coef\/f\/icients in $\mathbb{Z}[\alpha,\rho,\kappa]$. The remainder of the proof then follows from specializing the formulas in equation~\eqref{SPCD2} to $\beta=\Lambda_6=0$.
\end{proof}

The polynomial in equation~\eqref{SexticPolynomial2} is a twist of the polynomial given by
\begin{gather*}
 \hat{F}(x) = x^6 + a x^4 + b x^2 +1.
\end{gather*}
The curve $y^2= \hat{F}(x)$ has extra involutions, i.e., it has automorphisms other than the hyperelliptic involution, for appropriate values of $a$, $b$ (the discriminant is nonzero). In \cite{SV} for curves with automorphism the dihedral invariants
\begin{gather*}
u = a b, \qquad v =a^3+b^3,
\end{gather*}
were def\/ined which give a birational parametrization of this locus $\L_2$ which is a~two-dimensional subvariety of $\M_2$.
We have the following:
\begin{cor}
For the genus-two curve $\mathcal{C}\colon y^2=F(x)$ given by equation~\eqref{SexticPolynomial2} with $\chi_{35}=0$ we obtain the dihedral invariants
\begin{gather*}
 u  = \frac{\big(d^{(1)}_2 + d^{(2)}_2\big)\big(d^{(1)}_2 - d^{(2)}_2\big)}
 {\big(d^{(1)}_0 + d^{(2)}_0\big)\big(d^{(1)}_0 - d^{(2)}_0\big)},\\
 v = \frac{\big(d^{(1)}_2 - d^{(2)}_2\big)^3}{\big(d^{(1)}_0 + d^{(2)}_0\big)\big(d^{(1)}_0 - d^{(2)}_0\big)^2}
 + \dfrac{\big(d^{(1)}_2 + d^{(2)}_2\big)^3}{\big(d^{(1)}_0 + d^{(2)}_0\big)^2\big(d^{(1)}_0 - d^{(2)}_0\big)},
\end{gather*}
and the Igusa invariants $[J_2 : J_4 : J_6 : J_{10}]$ given by {\rm \cite[equation~(16)]{SV}}.
\end{cor}

\pdfbookmark[1]{References}{ref}
\LastPageEnding


\begin{thebibliography}{99}
\footnotesize\itemsep=0pt

\bibitem{data}
Beshaj L., Hidalgo R., Malmendier A., Kruk S., Quispe S., Shaska T., Rational points on the moduli space of genus two, in Algebraic Curves and their Fibrations in Mathematical Physics and Arithmetic Geometry,
 \textit{Contemporary Math.}, Vol.~703, Amer. Math. Soc., Providence, RI, 2018,
 87--120.

\bibitem{MR1505464}
Bolza O., On binary sextics with linear transformations into themselves,
 \href{https://doi.org/10.2307/2369402}{\textit{Amer.~J. Math.}} \textbf{10} (1887), 47--70.

\bibitem{MR3540958}
Booker A.R., Sijsling J., Sutherland A.V., Voight J., Yasaki D., A database of
 genus-2 curves over the rational numbers, \href{https://doi.org/10.1112/S146115701600019X}{\textit{LMS~J. Comput. Math.}}
 \textbf{19} (2016), suppl.~A, 235--254, \href{https://arxiv.org/abs/1602.03715}{arXiv:1602.03715}.

\bibitem{MR3349314}
Br\"oker R., Howe E.W., Lauter K.E., Stevenhagen P., Genus-2 curves and
 {J}acobians with a given number of points, \href{https://doi.org/10.1112/S1461157014000461}{\textit{LMS~J. Comput. Math.}}
 \textbf{18} (2015), 170--197, \href{https://arxiv.org/abs/1403.6911}{arXiv:1403.6911}.

\bibitem{MR0485106}
Clebsch A., Gordan P., Theorie der {A}belschen {F}unctionen, \textit{Thesaurus
 Mathematicae}, Vol.~7, Physica-Verlag, W\"urzburg, 1967.

\bibitem{MR871067}
Freitag E., Siegelsche {M}odulfunktionen, \href{https://doi.org/10.1007/978-3-642-68649-8}{\textit{Grundlehren der
 Mathematischen Wissenschaften}}, Vol.~254, Springer-Verlag, Berlin, 1983.

\bibitem{MR2899960}
Goren E.Z., Lauter K.E., Genus 2 curves with complex multiplication,
 \href{https://doi.org/10.1093/imrn/rnr052}{\textit{Int. Math. Res. Not.}} \textbf{2012} (2012), 1068--1142,
 \href{https://arxiv.org/abs/1003.4759}{arXiv:1003.4759}.

\bibitem{MR1438983}
Gritsenko V.A., Nikulin V.V., Igusa modular forms and ``the simplest''
 {L}orentzian {K}ac--{M}oody algebras, \href{https://doi.org/10.1070/SM1996v187n11ABEH000171}{\textit{Sb. Math.}} \textbf{187} (1996),
 1601--1641.

\bibitem{mod-curves}
Harris J., Morrison I., Moduli of curves, \href{https://doi.org/10.1007/b98867}{\textit{Graduate Texts in
 Mathematics}}, Vol.~187, Springer-Verlag, New York, 1998.

\bibitem{MR0114819}
Igusa J.I., Arithmetic variety of moduli for genus two, \href{https://doi.org/10.2307/1970233}{\textit{Ann. of Math.}}
 \textbf{72} (1960), 612--649.

\bibitem{MR0141643}
Igusa J.I., On {S}iegel modular forms of genus two, \href{https://doi.org/10.2307/2372812}{\textit{Amer.~J. Math.}}
 \textbf{84} (1962), 175--200.

\bibitem{MR0229643}
Igusa J.I., Modular forms and projective invariants, \href{https://doi.org/10.2307/2373243}{\textit{Amer.~J. Math.}}
 \textbf{89} (1967), 817--855.

\bibitem{MR527830}
Igusa J.I., On the ring of modular forms of degree two over {${\bf Z}$},
 \href{https://doi.org/10.2307/2373943}{\textit{Amer.~J. Math.}} \textbf{101} (1979), 149--183.

\bibitem{vishy}
Krishnamoorthy V., Shaska T., V\"olklein H., Invariants of binary forms, in
 Progress in {G}alois Theory, \href{https://doi.org/10.1007/0-387-23534-5_6}{\textit{Dev. Math.}}, Vol.~12, Springer, New
 York, 2005, 101--122, \href{https://arxiv.org/abs/1209.0446}{arXiv:1209.0446}.

\bibitem{MR3435723}
Lauter K., Naehrig M., Yang T., Hilbert theta series and invariants of genus~2
 curves, \href{https://doi.org/10.1016/j.jnt.2015.02.020}{\textit{J.~Number Theory}} \textbf{161} (2016), 146--174.

\bibitem{MR3366121}
Malmendier A., Morrison D.R., K3 surfaces, modular forms, and non-geometric
 heterotic compactif\/ications, \href{https://doi.org/10.1007/s11005-015-0773-y}{\textit{Lett. Math. Phys.}} \textbf{105} (2015),
 1085--1118, \href{https://arxiv.org/abs/1406.4873}{arXiv:1406.4873}.

\bibitem{satake}
Malmendier A., Shaska T., The {S}atake sextic in {F}-theory, \href{https://doi.org/10.1016/j.geomphys.2017.06.010}{\textit{J.~Geom.
 Phys.}} \textbf{120} (2017), 290--305, \href{https://arxiv.org/abs/1609.04341}{arXiv:1609.04341}.

\bibitem{Me}
Mestre J.F., Construction de courbes de genre {$2$} \`a partir de leurs
 modules, in Ef\/fective Methods in Algebraic Geometry ({C}astiglioncello,
 1990), \href{https://doi.org/10.1007/978-1-4612-0441-1_21}{\textit{Progr. Math.}}, Vol.~94, Birkh\"auser Boston, Boston, MA, 1991,
 313--334.

\bibitem{SV}
Shaska T., V\"olklein H., Elliptic subf\/ields and automorphisms of genus~2
 function f\/ields, in Algebra, Arithmetic and Geometry with Applications
 ({W}est {L}afayette, {IN}, 2000), \href{https://doi.org/10.1007/978-3-642-18487-1_42}{Springer}, Berlin, 2004, 703--723,
 \href{https://arxiv.org/abs/math.AG/0107142}{math.AG/0107142}.

\end{thebibliography}
\end{document}